\magnification=\magstep1
\hsize=16truecm
 
\input amstex
\TagsOnRight
\parindent=20pt
\parskip=2.5pt plus 1pt
\define\({\left(}
\define\){\right)}
\define\[{\left[}
\define\]{\right]}
\define\e{\varepsilon}

\define\supp {\sup\limits}

\define\summ{\sum\limits}
\define\prodd{\prod\limits}

\font\kisbetu=cmr8

\centerline{\bf A MULTIVARIATE GENERALIZATION OF}
\centerline{\bf HOEFFDING'S INEQUALITY}
\smallskip
\centerline{\it P\'eter Major}
\centerline{Alfr\'ed R\'enyi Mathematical Institute of the Hungarian
Academy of Sciences}
\centerline{Budapest, P.O.B. 127 H--1364, Hungary, e-mail:
major\@renyi.hu}
\medskip
 
{\narrower \noindent {\it Summary:}\/ We prove a multivariate
version of Hoeffding's inequality about the distribution of
homogeneous polynomials of Rademacher functions. The proof is
based on such an estimate about the moments of homogeneous
polynomials of Rademacher functions which can be considered as
an improvement of Borell's inequality in a most important
special case. \par}
 
\beginsection 1. Introduction. Formulation of the main results.
 
Hoeffding's inequality states the following result. (see e. g. [2],
Proposition 1.3.5.)
 \medskip\noindent
{\bf Theorem A. (Hoeffding's inequality).} {\it Let $\e_1,\dots,\e_n$
be independent random variables, $P(\e_j=1)=P(\e_j=-1)=\frac12$, $1\le
j\le n$, and let $a_1,\dots,a_n$ be arbitrary real numbers. Put
$Z=\summ_{j=1}^na_j\e_j$ and $V^2=\summ_{j=1}^na^2_j$. Then
$$
P(Z>u)\le\exp\left\{-\frac{u^2}{2V^2 }\right\}\quad
\text{for all }u>0. \tag1.1
$$
}\medskip
 
In the study of $U$-statistics we need a multivariate version
of this result. The goal of this paper is to present such an
inequality. To formulate it first we have to introduce some
notations.
 
Let us fix a positive integer~$k$ and some real numbers
$a(j_1,\dots,j_k)$ for all sets of arguments $\{j_1,\dots,j_k\}$
such that $1\le j_l\le n$, $1\le l\le k$, and $j_l\neq j_{l'}$ if
$l\neq l'$, in such a way that the numbers $a(j_1,\dots,j_k)$ are
symmetric functions of their arguments, i.e.
$a(j_1,\dots,j_k)=a(j_{\pi(1)},\dots,j_{\pi(k)})$ for all
permutations $\pi\in \Pi_k$ of the set $\{1,\dots,k\}$.
 
Let us define with the help of the above real numbers and a
sequence of independent random variables $\e_1,\dots,\e_n$,
$P(\e_j=1)=P(\e_j=-1)=\frac12$, $1\le j\le n$, the random variable
$$
Z=\sum\Sb (j_1,\dots, j_k)\: 1\le j_l\le n \text{ for all } 1\le l\le
k\\ j_l\neq j_{l'} \text{ if }l\neq l' \endSb a(j_1,\dots, j_k)
\e_{j_1}\cdots \e_{j_k} \tag1.2
$$
and the number
$$
V^2=\sum\Sb (j_1,\dots, j_k)\: 1\le j_l\le n \text{ for all } 1\le l\le
k\\ j_l\neq j_{l'} \text{ if }l\neq l' \endSb a^2(j_1,\dots, j_k).
\tag1.3
$$
Now we formulate the following result.
\medskip\noindent
{\bf Theorem 1. (The multivariate version of Hoeffding's inequality).}
{\it The random variable $Z$ defined in formula (1.2) satisfies the
inequality
$$
P(|Z|>u)\le A \exp\left\{-\frac12\(\frac uV\)^{2/k}\right\}
\quad\text{for all }u\ge 0 \tag1.4
$$
with the constant $V$ defined in (1.3) and some constants $A>0$
depending only on the parameter $k$ in the expression $Z$.}
\medskip
 
Let us remark that the condition that the coefficients
$a(j_1,\dots,j_k)$ are symmetric functions of their variables does
not mean a real restriction, since by replacing all coefficients
$a(j_1,\dots,j_k)$ by
$a_{\text{Sym}}(j_1,\dots,j_k)=\frac1{k!}\summ_{\pi\in\Pi_k}
a(j_{\pi(1)},\dots,j_{\pi(k)})$ in formula (1.2), where $\Pi_k$
denotes the set of all permutations of the set $\{1,\dots,k\}$ we
do not change the random variable $Z$. The identities $EZ=0$,
$EZ^2=k!V^2$ hold. A comparison of Theorem~A and Theorem~1 shows
that Theorem 1 yields a slightly weaker estimate in the special
case $k=1$ because of the pre-exponential coefficient $A$ in the
estimate (1.4). But the expressions in the exponent agree in
formula (1.1) and in formula (1.4) in the special case $k=1$.
 
Moreover, estimate (1.4), disregarding the pre-exponential
coefficient $A$ in it, is sharp for all  parameters $k\ge1$. To
see this let us consider the random variable $Z=Z_n$ defined in
(1.2) with the special choice
$$
a(j_1,\dots,j_k)=a_n(j_1,\dots,j_k)=\frac
V{\sqrt{n(n-1)\cdots(n-k+1)}}.
$$
It is known (see e.g.~[3]) that the random  variables $Z_n$
converge, as $n\to\infty$, in distribution to a random variable
which can be expressed by means of a $k$-fold Wiener--It\^o
integral. Moreover, it can be expressed in a more explicit form
as the distribution of $V\cdot H_k(\eta)$, where $\eta$ is a
random variable with standard normal distribution, and
$H_k(\cdot)$ is the $k$-th Hermite polynomial with leading
coefficient~1. Beside this, the tail behaviour of $H_k(\eta)$ is
similar to that of $\eta^k$ in a neighbourhood of the infinity.
Hence the above example shows that if we have no additional
restriction about the coefficients $a(j_1,\dots,j_k)$ of the random
variable $Z$, then the estimate (1.4) is essentially sharp. We
cannot write a better expression in the exponent of its right-hand
side. This problem is discussed in more detail in a more general
context in Example~2 of paper~[5].
 
Theorem 1 can be interpreted in such a way that the distribution of
$Z$ satisfies an inequality similar to the distribution of $V\eta^k$,
where $\eta$ is a standard normal random variable. We shall prove
it as a relatively simple consequence of the following result,
which formulates a similar statement about the moments of the
random variable~$Z$.
\medskip\noindent
{\bf Theorem 2.} {\it The random variable $Z$ defined in formula (1.2)
satisfies the inequality
$$
EZ^{2M}\le 1\cdot3\cdot5\cdots(2kM-1)V^{2M}\quad\text{for all }
M=1,2,\dots \tag1.5
$$
with the constant $V$ defined in formula (1.3).}
\medskip
We shall prove Theorem 2 with the help of two lemmas. Before their
formulation we introduce the following notation:
 
$$
\bar Z=\sum\Sb (j_1,\dots, j_k)\: 1\le j_l\le n \text{ for all } 1\le
l\le k\\ j_l\neq j_{l'} \text{ if }l\neq l' \endSb |a(j_1,\dots, j_k)|
\eta_{j_1}\cdots \eta_{j_k}, \tag1.6
$$
where $\eta_1,\dots,\eta_n$ are iid. random variables with standard
normal distribution, and the numbers $a(j_1,\dots,j_k)$ agree with
those in formula (1.2). Now we state
\medskip\noindent
{\bf Lemma 1.}
$$
EZ^{2M}\le E\bar Z^{2M}\quad\text{for all }M=1,2,\dots, \tag1.7
$$
\medskip\noindent
and
\medskip\noindent
{\bf Lemma 2.} {\it The random variable $\bar Z$ defined in formula
(1.6) satisfies the inequality
$$
E\bar Z^{2M}\le 1\cdot3\cdot5\cdots(2kM-1)V^{2M}\quad\text{for all }
M=1,2,\dots \tag1.8
$$
with the constant $V$ defined in formula (1.3).}
\medskip
 
Theorem 2 states an estimate about the moments of homogeneous
polynomials of the independent random variables $\e_1,\dots,\e_n$
which are sometimes called Rademacher functions in the literature.
We finish the Introduction by recalling Borell's inequality (see
e.g.~[1]) which gives a similar estimate. The proof of the
results will be given in Section 2. Then we compare Borell's
inequality with our results and make some comments in Section~3.
\medskip\noindent
{\bf Theorem B. (Borell's inequality).}
{\it The moments of the random variable $Z$ defined in formula (1.2)
satisfy the inequality
$$
E|Z|^p\le\(\frac{p-1}{q-1}\)^{kp/2} \(E|Z|^q\)^{p/q}\quad \text{ if }
\quad 1<q\le p<\infty. \tag1.9
$$
}
 
\beginsection 2. Proof of the results.
 
{\it Proof of Lemma 1.}\/ We can write, by carrying out the
multiplications in the expressions $EZ^{2M}$ and $E\bar Z^{2M}$,
by exploiting the additive and multiplicative properties of the
expectation for sums and products of independent random variables
together with the identities $E\e_j^{2k+1}=0$ and $E\eta_j^{2k+1}=0$
for all $k=0,1,\dots$  that
$$
EZ^{2M}=\sum\Sb j_1,\dots, j_l,\, m_1,\dots, m_l,\; 1\le j_s\le n\\
j_s\ge1,\; 1\le s\le l,\;\; m_1+\dots+m_l=M\endSb
A(j_1,\dots,j_l,m_1,\dots,m_l)E\e_{j_1}^{2m_1}\cdots E\e_{j_l}^{2m_l}
\tag2.1
$$
and
$$
E\bar Z^{2M}=\sum\Sb j_1,\dots, j_l,\, m_1,\dots, m_l,\; 1\le j_s\le n\\
j_s\ge1,\; 1\le s\le l,\;\; m_1+\dots+m_l=M\endSb
B(j_1,\dots,j_l,m_1,\dots,m_l)E\eta_{j_1}^{2m_1}\cdots
E\eta_{j_l}^{2m_l} \tag2.2
$$
with some coefficients $A(j_1,\dots,j_l,m_1,\dots,m_l)$ and
$B(j_1,\dots,j_l,m_1,\dots,m_l)$ such that
$$|
A(j_1,\dots,j_l,m_1,\dots,m_l)|\le
B(j_1,\dots,j_l,m_1,\dots,m_l). \tag2.3
$$
We could express the coefficients  $A(\cdot,\cdot,\cdot)$ and
$B(\cdot,\cdot,\cdot)$ in an explicit form, but we do not have to do
this. What is important for us is that $A(\cdot,\cdot,\cdot)$ can be
expressed as the sum of certain terms, and $B(\cdot,\cdot,\cdot)$ as
the sum of the absolute value of the same terms, hence relation (2.3)
holds. Since $E\e_j^{2m}\le E\eta_j^{2m}$ for all parameters $j$ and
$m$ formulas (2.1), (2.2) and (2.3) imply Lemma~1.
\medskip\noindent
{\it Proof of Lemma~2.} Let us consider a white noise $W(\cdot)$
on the unit interval $[0,1]$, i.e.\ let us take a set of Gaussian
random variables $W(A)$ indexed by the measurable sets $A\subset [0,1]$
such that $EW(A)=0$, $EW(A)W(B)=\lambda(A\cap B)$ with the Lebesgue
measure $\lambda$ for all measurable subsets of the interval $[0,1]$.
(We also need the relation $W(A\cup B)=W(A)+W(B)$ with probability 1
if $A\cap B=\emptyset$, but this relation is the consequence of the
previous ones. Indeed, they yield that $E(W(A\cup B)-W(A)-W(B))^2=0$
if $A\cap B=\emptyset$, and this implies the desired identity.)
Let us introduce the random variables $\eta_j=n^{1/2}W\(\[\frac
{j-1}n,\frac jn\)\)$, $1\le j\le n$, together with the function
$f(t_1,\dots,t_k)$, with arguments $0\le t_s<1$ for all indices
$1\le s\le k$, defined as
$$
f(t_1,\dots,t_k)=\left\{
\aligned
&n^{k/2}|a(j_1,\dots,j_k)| \quad\text{if }
t_s\in\[\frac {j_s-1}n,\frac {j_s}n\),\text{ and }j_s\neq j_{s'}
\text{ if } s\neq s', \\
&\hskip 7.5truecm  1\le j_s\le n,\;  1\le s\le k\\
&0 \quad\text{if }
t_s\in\[\frac {j_s-1}n,\frac {j_s}n\),\text{ and }j_s= j_{s'}
\text{ for some } s\neq s',\\
&\hskip 7.5truecm  1\le j_s\le n,\;  1\le s\le k
\endaligned \right. \tag2.4
$$
 
Observe that the above defined random variables $\eta_1,\dots,\eta_n$
are independent with standard normal distribution, hence we may
assume that they appear in the definition of the random variable
$\bar Z$ in formula (1.6). With such a choice we can represent
$\bar Z$ in the form of a $k$-fold Wiener--It\^o integral (introduced
e.g. in [4])
$$
\bar Z=\int f(t_1,\dots,t_k)W(\,dt_1)\dots W(\,dt_k)
$$
of the (elementary) function $f$ defined in formula (2.4) with
respect to white noise $W(t)$ we have introduced. Beside this,
the identity
$$
\int f^2(t_1,\dots,t_k)\,dt_1\dots\,dt_k=V^2
$$
also holds with the number~$V$ defined in formula (1.3). Hence to
complete the proof of Lemma~2 it is enough to show that if a
function~$f$ of $k$ variables and a $\sigma$-finite measure $\mu$
on some measurable space $(X,\Cal X)$ satisfy the inequality
$$
\int f^2(x_1,\dots,x_k)\mu(\,dx_1)\dots\mu(\,dx_k)=\sigma^2<\infty
$$
with some $\sigma^2>0$, then the moments of the $k$-fold
Wiener--It\^o integral (defined e.g. in~[4])
$$
J_{\mu,k}(f)=\frac1{k!}\int
f(x_1,\dots,x_k)\mu_W(\,dx_1)\dots\mu_W(\,dx_k)
$$
of the function $f$ with respect to a white-noise $\mu_W$ with
counting measure $\mu$ satisfy the inequality
$E\(k!J_{\mu,k}(f)\)^{2M}\le 1\cdot3\cdots(2kM-1)\sigma^{2M}$ for
all $M=1,2,\dots$. But this result (which can be got relatively
simply from the diagram formula for the product of Wiener--It\^o
integrals) is proven in Proposition~A of paper~[5], hence here I
omit the proof.
\plainfootnote{$^1$}{\kisbetu For the sake of completeness I put
the proof of this result together with some definitions needed to
understand it to an Appendix of this paper, but probably it will
not belong to the final version of this work.}
\medskip
Theorem 2 is a straightforward consequence of Lemmas~1 and~2.
Hence it remained to prove Theorem~1 with the help of Theorem~2.
\medskip\noindent
{\it Proof of Theorem~1.} By the Stirling formula we get from the
estimate of Theorem 2 that
$$
EZ^{2M}\le \frac{(2kM)!}{2^{kM}(kM)!} V^{2M}\le A
\(\frac2e\)^{kM}(kM)^{kM}V^{2M}
$$
for any $A\ge\sqrt2$ if $M\ge M_0(A)$. Hence we can write by the
Markov inequality that
$$
P(Z>u)\le\frac{EZ^{2M}}{u^{2M}}\le A\(\frac{2kM}e\(\frac
Vu\)^{2/k}\)^{kM}     \tag2.5
$$
for all $A>\sqrt 2$ if $M\ge M_0(A)$. Put $k\bar M=k\bar
M(u)=\frac12\(\frac uV\)^{2/k}$, and $M=M(u)=[\bar M]$, where
$[x]$ denotes the integer part of the number~$x$. Let us choose
a number $u_0$ by the identity $M(u_0)=M_0(A)$. Formula (2.5) can
be applied with $M=M(u)$ for $u\ge u_0$, and it yields that
$$
P(Z>u)\le Ae^{-kM}\le Ae^ke^{-k\bar M}=Ae^k\exp\left\{-\frac12
\(\frac uV\)^{2/k}\right\}\ \quad\text{if } u\ge u_0. \tag2.6
$$
Formula (2.6) means that relation (1.2) holds for $u\ge u_0$ if the
constant $A$ is replaced by $Ae^k$ in it. By choosing the constant
$A$ sufficiently large we can guarantee that relation (1.2) holds
for all $u\ge0$.
 
\beginsection 3. A discussion about the results.
 
Let us look what kind of estimate yields Borell's inequality for
the expression $Z$ defined in (1.2). It is natural to apply it with
the choice $q=2$. Since $EZ^2=k!V^2$, Borell's inequality yields
with  such a choice the estimate $E|Z|^{2p}\le (2p-1)^{kp}\(k!V^2\)^p$
for all real numbers $p\ge1$. Let us compare this inequality for the
moments $EZ^{2M}$ with large integers $M$ with the estimate of
Theorem~2. If we disregard some constant factors not depending on
$M$ we get that this estimate is of order $(2M)^{kM}V^{2M}\cdot
(k!)^M$, while Theorem~2 yields an estimate of order
$(2M)^{kM}V^{2M}\cdot\(\frac ke\)^{kM}$. It can be seen that
$k!>\(\frac ke\)^k$ for all $k\ge1$. This means that Borell's
inequality shows that $E Z^{2M}\le C^M(kM)^{kM}V^{2M}$ for large
$M$ with a universal constant $C$ depending only on the parameter
$k$ in formula~(1.2), but it does not give the optimal choice for
the parameter~$C$. As a consequence, it implies a weakened version
$P(|Z|>u)\le A \exp\left\{-B\(\frac uV\)^{2/k}\right\}$
of the inequality of Theorem~1 with some universal constants $A$
and $B$, but it cannot yield the optimal choice for the number $B$.
In short, Theorem~2 is weaker than Borell's  inequality in that
respect that it compares only the second and $2M$-th moment of the
random variable~$Z$, but it yields a sharper bound. Hence it can be
more useful in certain applications.
 
Let us finally remark that actually we have proved a sharper result
than Theorems~1 and~2. In those results we have defined the random
variable $Z$ with the help of independent random variables $\e_j$
with distribution $P(\e_j=1)=P(\e_j=-1)=\frac12$. But the proof of
Theorems~1 and~2 also works without any change in the case of random
variables with other distributions. Let us formulate this result.
First I introduce the following notion.
\medskip\noindent
{\bf Definition of sub-Gaussian distributions.} {\it Let us call a
random variable $\xi$ or its distribution sub-Gaussian, if its moments
satisfy the relations $E\xi^{2M-1}=0$ and $E\xi^{2M}\le E\eta^{2M}$
for all $M=1,2,\dots$, where $\eta$ is a random variable with standard
normal distribution.}
\medskip
It is clear that a random variable with distribution
$P(\e=1)=P(\e=-1)=\frac12$ is sub-Gaussian. Because of some
symmetrization arguments applied in probability theory this seems to
be the most important example of sub-Gaussian random variables, but
the following result holds for all of them.
\medskip\noindent
{\bf Theorem 3.} {\it  Let $\e_1,\dots,\e_n$ be independent
sub-Gaussian random variables (with possibly different
distributions). Let us define the random variable $Z$ by formula
(1.2) by the replacement of the original random variables
$\e_1,\dots,\e_n$ with these new random variables
$\e_1,\dots,\e_n$. This new random variable $Z$ also satisfies
the estimate (1.4) of Theorem~1 and the estimate (1.5) of
Theorem~2.}
 
\medskip
Theorem 3 means that the distribution and moments of homogeneous
polynomials of independent sub-Gaussian random variables satisfy
such estimates as the distribution and moments of homogeneous
polynomials of Gaussian random variables. Here the sub-Gaussian
property plays a most essential role. In the case of homogeneous
polynomials of independent, but not necessarily sub-Gaussian
random variables the situation is much more complex. But this
problem will not be discussed here.
 
\beginsection Appendix
 
To prove the inequality formulated at the end of Lemma~2 we need a
result which expresses the expected value of the product of multiple
Wiener--It\^o integrals in an appropriate way. To formulate this
result which is the simple consequence of a basic result of the
theory of Wiener--It\^o integrals, the so-called diagram formula,
first I have to introduce some notations. Let me recall that given a
$\sigma$-finite measure $\mu$ on some measurable space $(X,\Cal X)$
we call a white noise with counting measure $\mu$ such a Gaussian
field $\mu_W(A)$, $A\in \Cal X$, indexed by the measurable sets of
$X$ which satisfies the relations $E\mu_W(A)=0$ and
$E\mu_W(A)\mu_W(B)=\mu(A\cap B)$ for all $A,B\in\Cal X$.
 
Let us have a $\sigma$-finite measure $\mu$ together with a white
noise $\mu_W$ with counting measure $\mu$ on $(X,\Cal X)$. Let us
consider $L$ real valued functions $f_l(x_1,\dots,x_{k_l})$ on
$(X^{k_l},\Cal X^{k_l})$ such that
$\int f^2_l(x_1,\dots,x_{k_l})\mu(\,dx_1)\dots\mu(\,dx_{k_l})<\infty$,
$1\le l\le L$. Let us consider the Wiener--It\^o integrals
$k_l!J_{\mu,k_l}(f_l)=\int f_l(x_1,\dots,x_{k_l})
\mu_W(\,dx_1)\dots\mu_W(\,dx_{k_l})$, $1\le l\le L$,
and let us describe
how the expected value $E\(\prodd_{l=1}^L k_l!J_{\mu,k_l}(f_l)\)$
can be calculated by means of the diagram formula.
 
For this goal let us introduce the following notations. Put
$$
F(x_{(l,j)}, 1\le l\le L,\,1\le j\le k_l)=\prod_{l=1}^L
f_l(x_{(l,1)},\dots,x_{(l,k_l)}), \tag A1
$$
and define a class of diagrams $\Gamma(k_1,\dots,k_L)$ in the
following way: Each diagram $\gamma\in\Gamma(k_1,\dots,k_L)$ is a
(complete, undirected) graph with vertices $(l,j)$, $1\le l\le L$,
$1\le j\le k_l$, and we shall call the set of vertices $(l,j)$
with a fixed index $l$ the $l$-th row of the graphs
$\gamma\in\Gamma(k_1,\dots,k_L)$. The graphs
$\gamma\in\Gamma(k_1,\dots,k_L)$ will have edges with the following
properties. Each edge connects vertices $(l,j)$ and $(l',j')$ from
different rows, i.e. $l\neq l'$ for the end-points of an edge. From
each vertex there starts exactly one edge. $\Gamma(k_1,\dots,k_L)$
contains all graphs $\gamma$ with such properties. If there is no
such graph, then $\Gamma(k_1,\dots,k_L)$ is empty.
 
Put $2N=\summ_{l=1}^Lk_l$. Then each $\gamma\in \Gamma(k_1,\dots,k_L)$
contains exactly $N$ edges. If an edge of the diagram $\gamma$
connects some vertex $(l,j)$ with some other vertex $(l',j')$,
$l'>l$, then we call $(l',j')$ the lower end-point of this edge,
and we denote the set of lower end-points of $\gamma$ by $\Cal
A_\gamma$ which has $N$ elements. Let us also introduce the
following function $\alpha_\gamma$ on the vertices of $\gamma$. Put
$\alpha_\gamma(l,j)=(l,j)$ if $(l,j)$ is the lower end-point of an
edge, and $\alpha_\gamma(l,j)=(l',j')$ if $(l,j)$ is connected with
the point $(l'j')$ by an edge of $\gamma$, and $(l',j')$ is the
lower end-point of this edge. Then we define the function
$$
\bar F_\gamma(x_{(l,j)}, \; (l,j)\in \Cal A_\gamma)=
F(x_{\alpha_\gamma(l,j)}, 1\le l\le L,\,1\le j\le k_l)
$$
with the function $F$ introduced in (A1), i.e. we replace the
argument $x_{(l,j)}$ by $x_{(l',j')}$ in the function $F$ if $(l,j)$
and $(l',j')$ are connected by an edge in $\gamma$, and $l'>l$.
Then we enumerate the lower end-points somehow, and define the function
$B_\gamma(r)$, $1\le r\le N$, such that $B_\gamma(r)$ is the $r$-th
lower end-point of the diagram $\gamma$. Write
$$
F_\gamma(x_1,\dots,x_N)=\bar F_\gamma(x_{B_\gamma(r)},\; 1\le r\le N)
$$
and
$$
F_\gamma=\int\cdots\int
F_\gamma(x_1,\dots,x_N)\mu(\,dx_1)\dots\mu(\,dx_N)
\quad\text{for all }  \gamma\in \Gamma(k_1,\dots,k_L).
$$
Now we  formulate the corollary of the diagram formula we need.
\medskip\noindent
{\bf Theorem B.} {\it With the above introduced notation
$$
E\(\prodd_{l=1}^L k_l!J_{\mu,k_l}(f_l)\)=\sum_{\gamma\in
\Gamma(k_1,\dots,k_L)} F_\gamma.
$$
(If $\Gamma(k_1,\dots,k_L)$ is empty, then the expected value of
the above product of random integrals equals zero.) Beside this
$$
F_\gamma^2\le \prod_{l=1}^L\int
f^2_l(x_1,\dots,x_{k_l})\mu(\,dx_1)\dots\mu(dx_{k_l}) \quad
\text{for all } \gamma\in\Gamma(k_1,\dots,k_L).
$$
}\medskip
Now we turn to the proof of the inequality
$$
E\(k!J_{\mu,k}(f)\)^{2M}\le 1\cdot 3\cdot 5\cdots (2kM-1)\(\int
f^2(x_1,\dots,x_k)\mu(\,dx_1)\dots\mu(\,dx_k)\)^M. \tag A2
$$
\medskip\noindent
{\it Proof of Relation (A2).}\/ Relation (A2) can be simply proved
with the help of Theorem~B if we apply it with $L=2M$ and the
functions
$f_l(x_1,\dots,x_{k_l})=f(x_1,\dots,x_k)$ for all $1\le l\le 2M$.
Then Theorem~B yields that
$$
E\(k!J_{\mu,k}(f)^{2M}\)\le\( \int
f^2(x_1,\dots,x_k)\mu(\,dx_1)\dots\mu(dx_k)\)^M|\Gamma_{2M}(k)|,
$$
where $|\Gamma_{2M}(k)|$ denotes the number of diagrams $\gamma$
in $\Gamma(\underbrace{k,\dots,k}_{2M\text{ times}})$. Thus
to complete the proof of relation~(A2) it is enough to show that
$|\Gamma_{2M}(k)|\le 1\cdot3\cdot5\cdots (2kM-1)$. But this can be
seen simply with the help of the following observation. Let
$\bar\Gamma_{2M}(k)$ denote the class of all graphs with vertices
$(l,j)$, $1\le l\le 2M$, $1\le j\le k$, such that from all vertices
$(l,j)$ exactly one edge starts, all edges connect different vertices,
but we also allow edges connecting vertices $(l,j)$ and $(l,j')$ with
the same first coordinate~$l$. Let $|\bar\Gamma_{2M}(k)|$ denote
the number of graphs in $\bar\Gamma_{2M}(k)$. Then clearly
$|\Gamma_{2M}(k)|\le |\bar\Gamma_{2M}(k)|$. On the other hand,
$|\bar\Gamma_{2M}(k)|=1\cdot3\cdot5\cdots(2kM-1)$. Indeed, let us
list the vertices of the graphs from $\bar\Gamma_{2M}(k)$ in an
arbitrary way. Then the first vertex can be paired with another
vertex in $2kM-1$ way, after this the first vertex from which no
edge starts can be paired with $2kM-3$ vertices from which no edge
starts. By following this procedure the next edge can be chosen
$2kM-5$ ways, and by continuing this calculation we get the desired
relation.
 
\beginsection References
 
\item{1.)} Borell, C. (1979) On the integrability of Banach space
valued Walsh polynomials. {\it S\'eminaire de Probabilit\'es XIII,
Lecture Notes in Math. 721}\/ 1--3. Springer, Berlin.
\item{2.)} Dudley, R. M. (1998)  {\it Uniform Central Limit
Theorems.}\/ Cambridge University Press, Cambridge U.K.
\item{3.)} Dynkin, E. B. and Mandelbaum, A. (1983) Symmetric
statistics, Poisson processes and multiple Wiener integrals. {\it
Annals of Statistics\/} {\bf 11}, 739--745
\item{4.)} Major, P. (1981) Multiple Wiener--It\^o integrals. {\it
Lecture Notes in Mathematics\/} {\bf 849}, Springer Verlag, Berlin
Heidelberg, New York,
\item{5.)} Major, P. (2004) On a multivariate version of Bernstein's
inequality. Submitted to {\it Ann. Probab.}
 
\vskip2truecm\noindent
Abbreviated title: Hoeffding's inequality

\bye